\documentclass [12pt]{article}

\usepackage{amsfonts}
\usepackage{amssymb}
\usepackage{epsfig, subfigure, verbatim}

\def\cstar{${}$\circledast${}$}

\renewcommand{\emph}[1]{{\it #1}}
\newcounter{segcount}[section]
\renewcommand{\thesegcount}{\thesection.\arabic{segcount}}

\newenvironment{segment}
{\refstepcounter{segcount}\vspace{5mm}
\noindent{\bf \thesegcount. }}
{}

\newenvironment{statementnumbered}[3][]
{\refstepcounter{segcount}\vspace{5mm}
\noindent{\bf\thesegcount. #2}#1{\bf.}\ {\sl #3}}
{\nolinebreak[4] \nopagebreak[4] $\hfill \square$}
\newenvironment{statement}[3][]
{\vspace{5mm}\noindent{\bf #2}#1{\bf.} {\sl #3}}
{\nolinebreak[4] \nopagebreak[4] $\hfill \square$}
\newenvironment{statementnoboxnumbered}[3][]
{\refstepcounter{segcount}\vspace{5mm}
\noindent{\bf\thesegcount. #2}#1{\bf.} {\sl #3}}
{}
\newenvironment{statementnobox}[3][]
{\vspace{5mm}\noindent{\bf #2}#1{\bf.} {\sl #3}}
{}
\newenvironment{definitionnumbered}[2][]
{\refstepcounter{segcount}\vspace{5mm}
\noindent{\bf\thesegcount. #2}#1{\bf.}}
{}

\newenvironment{definition}[2][]
{\vspace{5mm}\noindent{\bf #2}#1{\bf.}}
{}

\newenvironment{resultnumbered}[3][]
{\refstepcounter{segcount}\vspace{5mm}
\noindent{\bf\thesegcount. #2}#1{\bf.} {\sl #3}
\vskip5mm\noindent {\bf Proof: }}
{\nolinebreak[4] \nopagebreak[4] $\hfill \square$}

\newenvironment{result}[3][]
{\vspace{5mm}
\noindent{\bf#2}#1{\bf.} {\sl #3}
{\\ \bf Proof: }}
{\nolinebreak[4] \nopagebreak[4] $\hfill \square$}

\newenvironment{risultnumbered}[3][]
{\refstepcounter{segcount}\vspace{5mm}
\noindent{\bf\thesegcount. #2}#1{\bf.} {\sl #3}
{\nopagebreak[4] \noindent \bf Proof: }}
{$\hfill \square$}

\newenvironment{risult}[3][]
{\vspace{5mm}
\noindent{\bf#2}#1{\bf.} {\sl #3}
{\noindent \bf Proof: }}
{\nolinebreak[4] \nopagebreak[4] $\hfill \square$}

\newenvironment{risultnoboxnumbered}[3][]
{\refstepcounter{segcount}\vspace{5mm}
\noindent{\bf\thesegcount. #2}#1{\bf.} {\sl #3}
{\nopagebreak[4] \noindent \bf Proof: }}
{}

\newenvironment{risultnobox}[3][]
{\vspace{5mm}
\noindent{\bf#2}#1{\bf.} {\sl #3}
{\nopagebreak[4] \noindent \bf Proof: }}
{}

\newenvironment{resultnoboxnumbered}[3][]
{\refstepcounter{segcount}\vspace{5mm}
\noindent{\bf\thesegcount. #2}#1{\bf.} {\sl #3}
{\\ \bf Proof: }}
{}

\newenvironment{resultnobox}[3][]
{\vspace{5mm}\noindent{\bf #2}#1{\bf.} {\sl #3}
{\\ \bf Proof: }}
{}

\newcommand{\seg}{\begin{segment}}
\newcommand{\segend}{\end{segment}}

\newcommand{\stmtnum}{\begin{statementnumbered}}
\newcommand{\stmtnumend}{\end{statementnumbered}}

\newcommand{\stmt}{\begin{statement}}
\newcommand{\stmtend}{\end{statement}}

\newcommand{\stmtnoboxnum}{\begin{statementnoboxnumbered}}
\newcommand{\stmtnoboxnumend}{\end{statementnoboxnumbered}}

\newcommand{\stmtnobox}{\begin{statementnobox}}
\newcommand{\stmtnoboxend}{\end{staytementnobox}}

\newcommand{\defnnum}{\begin{definitionnumbered}}
\newcommand{\defnnumend}{\end{definitionnumbered}}

\newcommand{\defn}{\begin{definition}}
\newcommand{\defnend}{\end{definition}}

\newcommand{\resnum}{\begin{resultnumbered}}
\newcommand{\resnumend}{\end{resultnumbered}}

\newcommand{\res}{\begin{result}}
\newcommand{\resend}{\end{result}}

\newcommand{\risnum}{\begin{risultnumbered}}
\newcommand{\risnumend}{\end{risultnumbered}}

\newcommand{\ris}{\begin{risult}}
\newcommand{\risend}{\end{risult}}

\newcommand{\risnoboxnum}{\begin{risultnoboxnumbered}}
\newcommand{\risnoboxnumend}{\end{risultnoboxnumbered}}

\newcommand{\risnobox}{\begin{risultnobox}}
\newcommand{\risnoboxend}{\end{risultnobox}}

\newcommand{\resnoboxnum}{\begin{resultnoboxnumbered}}
\newcommand{\resnoboxnumend}{\end{resultnoboxnumbered}}

\newcommand{\resnobox}{\begin{resultnobox}}
\newcommand{\resnoboxend}{\end{resultnobox}}

\newcommand{\td}{\downarrow}
\newcommand{\tD}{\Downarrow}
\newcommand{\la}{\langle}
\newcommand{\ra}{\rangle}

\newcommand{\cA}{{\cal A}}
\newcommand{\cB}{{\cal B}}
\newcommand{\cC}{{\cal C}}

\newcommand{\cG}{{\cal G}}
\newcommand{\cH}{{\cal H}}

\newcommand{\cM}{{\cal M}}

\newcommand{\cP}{{\cal P}}
\newcommand{\cQ}{{\cal Q}}
\newcommand{\cR}{{\cal R}}
\newcommand{\cS}{{\cal S}}
\newcommand{\cT}{{\cal T}}
\newcommand{\cU}{{\cal U}}

\newcommand{\bS}{{\bf S}}

\def\spsubseteq{\subset\hspace{-2.1mm}\raise0.45ex\hbox{.}\hspace{1mm}}
\def\spsubsetneqq{\subsetneqq\hspace{-2.1mm}\raise1ex\hbox{.}\hspace{1mm}}

\title{Unique Factorisation of Additive Induced-Hereditary Properties
}

\author{
Alastair Farrugia%
\thanks{\noindent The first author's studies in Canada are fully funded by the
Canadian government through a Canadian Commonwealth Scholarship.  
The second author's research is financially supported by NSERC.
The results presented here form part of the first author's Ph.D. thesis, that he is writing under the supervision of the second author. 
} 
\\ {\it afarrugia@math.uwaterloo.ca}
\and
R. Bruce Richter
\\ {\it brichter@math.uwaterloo.ca}
\\ Dept. of Combinatorics \& Optimization
\\ University of Waterloo, Ontario, Canada, N2L 3G1
}

\begin{document}

\maketitle

\begin{abstract}
An additive hereditary graph property is a set of graphs, 
closed under isomorphism and under taking subgraphs and disjoint unions.
Let $\cP_1, \ldots, \cP_n$ be additive hereditary graph properties. A
graph $G$ has property $(\cP_1 \circ \cdots \circ \cP_n)$ if there is
a partition $(V_1, \ldots, V_n)$ of $V(G)$ into $n$ sets such that,
for all $i$, the 
induced subgraph $G[V_i]$ is in $\cP_i$. 
A property $\cP$ is reducible if there are properties $\cQ$, $\cR$
such that $\cP = \cQ \circ \cR$; otherwise it is
irreducible. Mih\'{o}k, Semani\v{s}in and Vasky [J. Graph Theory {\bf 33}
(2000), 44--53] gave a factorisation 
for any additive hereditary property $\cP$ into a given number $dc(\cP)$
of irreducible additive hereditary factors. 
Mih\'{o}k  [Discuss.\ Math.\ Graph Theory {\bf 20} (2000), 143--153] gave a similar factorisation for properties that are additive and induced-hereditary (closed under taking induced-subgraphs and disjoint unions). Their results left open the possiblity of different factorisations, maybe even with a different number of factors; we prove here that the given factorisations are, in fact, unique.
\end{abstract}

\section{Introduction\label{sec-intro}}
A {\em graph property\/} is an isomorphism-closed set of graphs. 
A graph $G$ {\em has property $\cP$\/} if $G \in \cP$. 
The {\em universal property\/}
${\cal U}$ is the set of all (finite, unlabelled, simple) graphs. A
property $\cP$ is {\em non-trivial\/} if $\emptyset \not= \cP \not= {\cal U}$. 

A property is {\em hereditary\/}, {\em induced-hereditary\/} or {\em
additive\/}  if
it is closed under taking subgraphs, induced-subgraphs or disjoint
unions, respectively. If $\cP$ is additive, and every %connected
component of a graph $X$ is in $\cP$, then $X$ is also in $\cP$.  

Let $\cP_1, \ldots, \cP_n$ be graph properties. A {\em $(\cP_1, \ldots,
\cP_n)$-partition\/} of a graph $G$ is a partition $(V_1, \ldots, V_n)$
of $V(G)$ into $n$ (possibly empty) sets such that, for all $i$, the
induced subgraph $G[V_i]$ is in $\cP_i$. The property $\cP = \cP_1
\circ \cdots \circ \cP_n$ is the set of all graphs
having a $(\cP_1, \ldots, \cP_n)$-partition. 
The $\cP_i$'s are {\em factors\/} or {\em divisors\/} of
$\cP$, while 
$\cP$ is the {\em product\/} of the $\cP_i$'s. 
It is easy to see that the product of
additive (or hereditary or induced-hereditary) 
properties is also additive (or hereditary or induced-hereditary).

%A property is {\em reducible\/} if it is the product of at least 
%two non-trivial properties; otherwise it is {\em irreducible\/}.
In this article, an additive (induced-)hereditary property is {\em reducible} 
if it is the product of at least two non-trivial 
\emph{additive (induced-)hereditary} properties; 
otherwise it is {\em irreducible\/}.
We show in~\cite{uf-her} that if an additive (induced-)hereditary
property is the product of \emph{any} two non-trivial properties, then it
is also the product of two additive (induced-)hereditary non-trivial properties.
So the concept of reducibility used here turns out to be the same 
as a more natural concept of reducibility; we point out,
however, that the proofs in~\cite{uf-her} depend on this article.

Mih\'{o}k, Semani\v{s}in and Vasky~\cite{uft-1} gave a factorisation
of an additive hereditary property $\cP$ into a given number $dc(\cP)$
of irreducible additive hereditary factors. This factorisation was
shown to be well-defined, but it was also claimed to be unique. 
The argument was that if $\cP = \cQ \circ \cR$, where by induction
$\cQ$ and $\cR$ each have a unique factorisation, then $\cP$ also has
a unique factorisation. However, there is still the possibility that
$\cP$ factors as $\cP_1 \circ \cdots \circ \cP_r$, where no subset of
the $\cP_i$'s has either $\cQ$ or $\cR$ as a product. 

As an analogy,
consider the ring $\{x+y \sqrt{5} \mid x,y \in \mathbb{Z}\}$. In this
integral domain $2$, $1+\sqrt{5}$ and $1-\sqrt{5}$ are all irreducible ---
they have no factorisation into two non-unit factors. In particular,
$2$ has a unique factorisation, but $4 = 2^2$ does not, because
we have $4 = (1+\sqrt{5})(-1+\sqrt{5})$.

We show in Theorem~\ref{uft-n} that similar 
anomalies do not occur with additive hereditary properties if the two
factorisations have \emph{exactly} $dc(\cP)$ factors. Could there
then be factorisations with different numbers of factors? This is not
an idle question, as Mih\'ok et al.\ showed in
Example 4.2 of the same paper that a certain
hereditary (but not additive) property $\cP_1 \circ \cP_2$ 
has another factorisation $\cQ_1 \circ \cQ_2 \circ \cQ_3$ 
where even the number of irreducible hereditary factors is different. 
One of the main contributions of this paper is Theorem~\ref{m=n},
where we prove that in any
factorisation into irreducible additive hereditary factors, 
the number of factors of $\cP$ must be \emph{exactly} $dc(\cP)$.  
%***P.\ Mih\'ok has recently informed us of a somewhat simpler proof of this
%***result, that will appear in a joint paper.
%\change-this?
%%P.\ Mih\'ok has recently suggested a simpler, though not shorter, 
%%way of proving 
%%Theorem~\ref{m=n}, that we will publish in a joint paper.

In~\cite{uft-2}, Mih\'{o}k gave a remarkably general 
%+%constuction of
construction of
uniquely partitionable graphs, and used this to produce a
factorisation for the wider class of 
properties that are additive and induced-hereditary.  This was claimed
to be unique using the same argument as in \cite{uft-1}.  We generalise his
construction, and our own results
(Theorems~\ref{m-resp-n},~\ref{ind-her-n} and~\ref{ind-her-m=n},
respectively) 
to prove that this factorisation is in fact unique.

We note that unique factorisation was settled completely in~\cite{uft-0} for a
significant class of additive hereditary properties, the proof
depending on the structure of those properties (and in the spirit of
the proof we give here). It is possible to use the structure of the
factorisation 
presented in~\cite{uft-1} to show that any factorisation with exactly
$dc(\cP)$ additive hereditary factors must be the one constructed in
that article (a similar proof is possible for the factorisation
of~\cite{uft-2}); the appeal of the proofs of uniqueness given here 
%for Theorems~\ref{uft-n} and~\ref{ind-her-n} 
is that they are independent
of the structure of the factors of $\cP$. Thus, they depend only on
the more elementary aspects of \cite{uft-1} and \cite{uft-2}.

In the next section we reproduce the essential concepts, definitions
and results  
adapted from~\cite{uft-1}; stating those results in a
stronger fashion here, and sometimes omitting their simple proofs. 
Our own techniques and proofs are presented in Section~\ref{sec-proof}.  
Sections~\ref{sec-ind-def} and~\ref{sec-ind-proof} are the induced-hereditary
analogues of these two sections.  
We end with some corollaries of unique factorisation 
% in Section~\ref{sec-consequences}, 
and a list of open questions.
% in Section~\ref{sec-open}.

A second paper~\cite{uf-her} contains related results on 
uniquely partitionable graphs, a characterisation of
induced-hereditary properties uniquely factorisable into 
arbitrary properties (not necessarily induced-hereditary).
A technical report~\cite{corr} contains the results of both
papers, and generalises them. More recently, the first 
author~\cite{farr-comp} used the results in~\cite{uft-2} 
and in this paper to show that it is NP-hard to recognise
reducible additive induced-hereditary properties, with 
the exception of the set of bipartite graphs.

\section{Definitions and results from~\cite{uft-1}\label{sec-def}}
In this section and the next we will actually prove unique factorisation 
for a class of properties strictly larger than the additive hereditary class.
We use $G \subseteq H$ to denote that $G$ is a subgraph of $H$. 
A {\em hereditary compositive} property is a hereditary property $\cP$ 
where, for any two graphs $G_1, G_2 \in \cP$, there is a graph $H \in \cP$
such that $G_i \subseteq H, i = 1,2$. It turns out that the proof of
unique factorisation for additive hereditary properties carries over 
to the hereditary compositive case without any change.
For our purposes, a hereditary compositive property is reducible if it is the
product of two non-trivial hereditary compositive properties; otherwise it is irreducible.

The unique factorisation result
for additive induced-hereditary properties includes as a special
case the result for additive hereditary properties
(Prop.~\ref{dc=dec}), but not the one 
for hereditary compositive properties.
An additive hereditary property is both additive 
induced-hereditary, and hereditary compositive.
However, for a fixed finite graph $S$, 
properties of the form $\cP_S := \{G \mid G \subseteq S\}$ are
hereditary compositive but not additive.
In~\cite{corr} we prove unique factorisation for a class that strictly
includes additive induced-hereditary properties, but still does not
contain properties of the form $\cP_S$.

In addition, the structures of the proofs of 
Theorems \ref{uft} and \ref{ind-her-uft} are similar.  
Having in mind the simpler proof for Theorem \ref{uft} 
before attempting the more difficult proof of 
Theorem~\ref{ind-her-uft} is very helpful. 

The smallest hereditary property that contains a set $\cG$ is 
denoted by $[\cG]$.  This is the  hereditary property 
{\em generated by $\cG$\/}, or that {\em $\cG$ generates\/}. 
$\cG$ {\em is a generating
set for\/}  $\cP$ if $[\cG] = \cP$. 
It is easily seen that

\[\begin{array}{lll}
[\cG] &=& \{G \mid \exists\, H \in \cG,\ G \subseteq H \}.
\end{array}\]

The {\em completeness\/} $c(\cP)$ of a hereditary property $\cP \not=
{\cal U}$ is 
$\max\{k: K_k \in \cP\}$\footnote{In some of the literature,
the convention is to 
define $c(\cP) := \max\{k: K_{k+1} \in \cP\}$.}, where $K_k$ is the
complete graph on $k$ vertices; clearly, $c(\cQ \circ \cR) = c(\cQ) +
c(\cR)$. Thus, any factorisation of a hereditary property $\cP$ has at most $c(\cP)$ 
non-trivial factors.

The {\em join\/} $G_1 + \cdots + G_n$ of $n$ graphs $G_1, \ldots, G_n$
consists of disjoint copies of the $G_i$'s, and all edges between
$V(G_i)$ and $V(G_j)$, for $i \not= j$. A graph $G$ is {\em
decomposable\/} if it is the join of
two graphs; otherwise, $G$ is
{\em indecomposable\/}. It is easy to see that $G$ is decomposable if and only
if its complement $\overline{G}$ is disconnected; $G$ is the join of
the complements of the components of $\overline{G}$, so every
decomposable graph can be expressed uniquely as the join of
indecomposable subgraphs, the {\em ind-parts\/} of $G$. The number of
ind-parts of $G$ is the {\em decomposability number\/} $dc(G)$ of $G$.

For a hereditary
property $\cP$, a graph $G$ is {\em $\cP$-strict\/} if $G \in \cP$ but
$G+K_1 \not\in \cP$. 
The set $\cM(\cP)$ of $\cP$-maximal graphs is defined as: 

%\smallskip
\begin{eqnarray*}
\cM(n,\cP) &:=& \{G \in \cP \mid |V(G)| = n \textrm{ and for all } e \in 
E(\overline{G}), G + e \not\in \cP \}; \\
\cM(\cP) &:=& {\displaystyle \bigcup_{n=1}^{\infty} \cM(n,\cP)}
\end{eqnarray*}

\smallskip
Since,  for $1 \leq n \leq c(\cP)$, $M(n,\cP) = \{K_n\}$, it is also
useful to define 
\[ \cM^*(\cP) := {\displaystyle \bigcup_{n=c(\cP)}^{\infty} \cM(n,\cP)}. \]

\vspace{-5mm}
\stmtnum[\ \cite{uft-1}]{Lemma\label{max-char}}
{Let $\cP_1, \ldots, \cP_m$ be hereditary properties of graphs, and
denote $\cP_1 \circ \cdots \circ \cP_m$ by $\cP$. A graph $G$ belongs
to $\cM(\cP)$ if and only if, for every $(\cP_1, \ldots,
\cP_m)$-partition $(V_1, \ldots, V_m)$ of $V(G)$, the following holds: 
$G[V_i] \in \cM(\cP_i)$ for $i = 1,\ldots, n$, and $G = G[V_1] +
\cdots + G[V_m]$. Moreover, if $G \in \cM^*(\cP)$, then it is
$\cP$-strict,  each $G[V_i]$ is $\cP_i$-strict, and is in
$\cM^*(\cP_i)$; in particular, each $G[V_i]$ is non-empty. 
}
\stmtnumend
\newline

It follows that if $\cP$ is reducible, then every graph in
$\cM^*(\cP)$ is decomposable.  
We note that the join of a $\cQ$-maximal graph $G$ and an
$\cR$-maximal graph $H$ need not be $(\cQ \circ \cR)$-maximal; for
example, take $G$ to be complete, $|V(G)| \leq c(\cQ) - 2$, and $H$
not complete. 

Clearly $[\cM^*(\cP)] = \cP$, but if $\cP$ is additive it is not the
unique generating set. 
If $\cG$ is a generating set for the hereditary property $\cP$, its 
{\em decomposability number\/} $dc(\cG)$ is 
$\min\{dc(G) \mid G \in
\cG\}$;
 the {\em decomposability number\/} of $\cP$ is $dc(\cP) :=
dc(\cM^*(\cP))$.  
A property with $dc(\cP) = 1$ is {\em indecomposable\/}; by
Lemma~\ref{max-char} such a property must be irreducible, and we shall
see that for hereditary \emph{compositive} properties the converse is also true.
The converse is not true for hereditary properties in general, as shown
in~\cite{uft-1}.

\stmtnum[\ \cite{uft-1}]{Lemma\label{gen-0}}
{Let $\cP$ be a hereditary property and let 
$G \in \cM^*(\cP)$, $H \in \cP$. If $G \subseteq H$ then $dc(H) \leq dc(G)$. If we have equality, with $G = G_1 + \cdots + G_n$ and $H = H_1 + \cdots + H_n$ being the respective expressions as joins of ind-parts, then 
we can relabel the ind-parts of $H$ so that each
$G_i$ is an induced subgraph of $H_i$.
} 
\stmtnumend
%\newline

\stmtnum[\ \cite{uft-1}]{Lemma\label{gen-1}}
{If $\cG$ generates the hereditary property $\cP$, then $dc(\cG) \leq
dc(\cM^*(\cP))$, with equality if $\cG \subseteq \cM^*(\cP)$. 
}
\stmtnumend
\newline

For $\cG\subseteq \cP$ and $H\in\cP$, let $\cG[H] := \{G \in
\cG \mid H \subseteq G \}$.

\stmtnum[\ \cite{uft-1}]{Lemma\label{gen-2}}
{Let $\cG$ generate the hereditary compositive property
$\cP$, and let $H$ be an arbitrary graph in $\cP$. Then $\cG[H]$ also
generates $\cP$.
}
\stmtnumend
\newline

For a generating set $\cG\subseteq \cM^*(\cP)$, let $\cG^{\td} := \{G
\in \cG \mid dc(G) = dc(\cP) \}$. 

\stmtnum[\ \cite{uft-1}]{Lemma\label{gen-3}}
{If $\cG \subseteq \cM^*(\cP)$ generates the hereditary compositive
property $\cP$, then so does $\cG^{\td}$. 
}
\stmtnumend
%\newline

\section{Unique factorisation for hereditary \\
compositive properties\label{sec-proof}} 

Our interpretation of \cite{uft-1} is that Mih\'ok et al.\ proved that every
hereditary compositive property $\cP$ has a factorisation into $dc(\cP)$
indecomposable factors.  Therefore, reducibility and decomposability
are the same thing.  Our purpose here is to show that every hereditary 
compositive property has at most one factorisation into indecomposable
hereditary compositive factors. We do so in the following two results.

\stmtnoboxnum{Theorem\label{m=n}}
{Let 
$\cP_1 \circ \cdots \circ \cP_m$ be a factorisation of the
hereditary compositive property $\cP$ into indecomposable 
hereditary compositive properties. Then $m = dc(\cP)$. 
}\stmtnoboxnumend

\stmtnoboxnum{Theorem\label{uft-n}}
{A hereditary compositive property $\cP$ can have only one factorisation
with exactly $dc(\cP)$ indecomposable hereditary compositive factors. 
}\stmtnoboxnumend

\medskip
The following result from \cite{uft-1} shows there is at least one
factorisation.

\stmtnum[\ \cite{uft-1}]{Theorem\label{factorn}}
{A hereditary compositive property has a factorisation into $dc(\cP)$
(necessarily indecomposable) hereditary compositive factors. 
}
\stmtnumend
%\newline

\medskip This in turn implies the following.

\stmtnum[\ \cite{uft-1}]{Corollary\label{irr-ind}}
{A hereditary compositive property is irreducible if and only if it is
indecomposable. 
}
\stmtnumend
%\newline

\medskip Putting this all together, we conclude:

\stmtnum{Hereditary Compositive Unique Factorisation Theorem\label{uft}}
{A hereditary compositive property has a unique factorisation
into irreducible hereditary compositive factors, and the number of
factors is exactly $dc(\cP)$. 
}
\stmtnumend
%newline

\medskip Our proofs 
depend heavily on the following construction of a 
generating set for $\cP$. 
Suppose $\cP_1 \circ \cdots \circ\cP_m$ is a factorisation of $\cP$ 
into  indecomposable hereditary compositive factors, and, for each $i$, we
are given a generating set $\cG_i \subseteq \cM^*(\cP_i)$ and a graph
$H_i \in \cP_i$. By Lemmas~\ref{gen-2} and~\ref{gen-3}, the set
$\cG_i^{\td}[H_i] := \{G \in \cG_i \mid H_i \subseteq G,\ dc(G) = 1
\}$ is also a generating set for $\cP_i$.  

We set\footnote
{Our notation extends easily to the join of any $m$ sets: $\cG_1 +
\cdots + \cG_m$; and to generating sets that contain several specified
subgraphs: $\cG[H_1, \ldots, H_r]$.} 
$\cG_1^{\td}[H_1] + \cdots + \cG_m^{\td}[H_m] := \{ G_1 + \cdots + G_m
\mid \forall\,i\ G_i \in \cG_i^{\td}[H_i] \}$.  This is clearly a
generating set for $\cP$, but need not consist of $\cP$-maximal graphs
(even if $m = dc(\cP)$).   However, we can add edges to
each graph $G_1+\cdots+G_m$ until we get (in all possible ways) 
a $\cP$-maximal graph $G'$.
Using $G \spsubseteq H$ to mean that $G$ is a spanning
subgraph of $H$, we can now describe the generating set we want:  

\begin{eqnarray*}
(\cG_1[H_1] + \cdots + \cG_m[H_m])^{\td} & := & \{G' \in \cM^*(\cP) \mid
dc(G') = dc(\cP), \textrm{ and } \\ 
&& \exists\, G \in \cG_1^{\td}[H_1] + \cdots + \cG_m^{\td}[H_m]\, ,\  G
\spsubseteq G'\}. 
\end{eqnarray*}

The following is immediate from the definition, and from Lemma~\ref{gen-3}.

\stmtnoboxnum{Lemma}
{Let $\cG = (\cG_1[H_1] + \cdots + \cG_m[H_m])^{\td}$.  Then:
\begin{enumerate} 
\item$\cG$ is a generating set for $\cP = \cP_1 \circ \cdots \circ 
\cP_m$; 
\item if $G \in \cG$, then $dc(G) = dc(\cP)$; and 
\item every $G\in \cG$ is spanned by  a join of $m$ indecomposable graphs, each
of which contains a different one of $H_1, \ldots, H_m$.
$\hfill \square$
\end{enumerate}
}
\stmtnoboxnumend
%\newline

Because we take $G' = G'_1 + \cdots + G'_{dc(\cP)}\in\cG$ to be a
\emph{spanned} supergraph of $G = G_1 + \cdots + G_m\in
\cG_1^{\td}[H_1] + \cdots + \cG_m^{\td}[H_m]$,  we
must have, for each $i$,
 $V(G_i) = V(\sum_{j\in J_i} G'_j)$, where $(J_1, J_2,
\ldots, J_m)$ is some partition of $\{1,2,\ldots,n\}$. That is, each
of the $m$ ind-parts of $G$ is a spanning subgraph of a join of
ind-parts from $G'$. We note that although $G_i \in \cP_i$, none of
the $G'_j$, $j\in J_i$, need be in $\cP_i$. In particular, 
the crucial observation that
Theorem~\ref{m=n} rests on is that, if $|J_i| > 1$, then $G_i \spsubsetneqq
\sum_{j\in J_i} G'_j$, and, since $G_i$ was $\cP_i$-maximal, $\sum_{j\in
J_i} G'_j$ is not in $\cP_i$. 

%\resnum{Theorem\label{uft-n}}
%{A hereditary compositive property $\cP$ can have only one factorisation
%with exactly $dc(\cP)$ indecomposable hereditary compositive factors. 
%}

We present first the proof of Theorem \ref{uft-n}, since it is simpler.

\medskip\noindent{\bf Proof of Theorem \ref{uft-n}:} 
Let $\cP_1 \circ \cdots \circ \cP_n = \cQ_1 \circ \cdots \circ \cQ_n$ 
be two factorisations of $\cP$ into $n = dc(\cP)$
indecomposable hereditary compositive factors. 

Label the $\cP_i$'s inductively, beginning with $i=n$, so that, for each
$i$, $\cP_i$ is inclusion-wise maximal among
$\cP_1,\cP_2,\dots,\cP_i$.  For each $i,j$ such that $i>j$, if
$\cP_i\setminus \cP_j\ne\emptyset$, then let $X_{i,j}\in
\cP_i\setminus \cP_j$; if $\cP_i\setminus \cP_j=\emptyset$, then
$\cP_i= \cP_j$ and we set $X_{i,j}$ to be the null graph.  For each
$i$, by compositivity there is an $H_{i,0}\in\cP_i$ that contains
all the $X_{i,j}$'s as subgraphs. The important point is
that if $\{L_1, L_2,\dots, L_n\}$ is an
unordered $(\cP_1,\dots,\cP_n)$-partition of some graph $G$ such that, for
each $i=1,2,\dots,n$, $H_{i,0}\subseteq G[L_i]$, then, by reverse 
induction on $i$ starting at $n$, $G[L_i]\in\cP_i$.

For each $i$, let $\cG_i = \{G_{i,0}, G_{i,1},
G_{i,2}, \ldots \}$ be a generating set for $\cP_i$. 
When graphs have a double
subscript, we will use the second number to denote which step of our
construction we are in. We start with $H_0 = H_{1,0} + \cdots +
H_{n,0}$.

For each $s \geq 0$, let $H_{s+1} \in (\cG_1[H_{1,s},
G_{1,s}] + \cdots + \cG_n[H_{n,s}, G_{n,s}])^{\td}$. Then $H_{s+1}$ has
an ind-part from each $\cG_i[H_{i,s}, G_{i,s}]$; we label the
ind-parts as $H_{1,s+1}$, \dots, $H_{n,s+1}$, so that,
 for each $i$, $H_{i,1} \subseteq H_{i,2}
\subseteq H_{i,3} \subseteq \cdots $

For $\cG_i[H_{i,s}, G_{i,s}]$ to be non-empty, we must have $H_{i,s}
\in \cP_i$. We know that the $H_{i,s+1}$'s give an unordered $\{\cP_1,
\ldots, \cP_n\}$-partition of $H_{s+1}$. From the earlier remark, for
$i=1,2,\dots,n$, $H_{i,s+1}\in \cP_i$.

By Lemma~\ref{max-char}, the ind-parts of $H_s$ form its $\{\cQ_1,
\ldots, \cQ_n\}$-partition, so there is some permutation $\varphi_s$
of $\{1,2,\ldots, n\}$ such that, for each $i$, $H_{i,s} \in
\cQ_{\varphi_s(i)}$. Since there are only finitely many permutations
of $\{1,2,\ldots,n\}$, there must be some permutation $\varphi$ that
appears infinitely often. Now whenever $\varphi_t = \varphi$, we have
$H_{i,1} \subseteq H_{i,2} \subseteq \cdots \subseteq H_{i,t} \in
\cQ_{\varphi(i)}$, so by heredity, for every $s\le t$, $H_{i,s}$ is in
$\cQ_{\varphi(i)}$.  Therefore,  we can take $\varphi_s = \varphi$ for all
$s$. By re-labelling the $\cQ_i$'s, we can assume $\varphi$ is the
identity permutation, so that $H_{i,s} \in \cQ_i$ for all $i$ and
$s$. 

For each $i$ and $s$, $G_{i,s-1} \subseteq H_{i,s}$, so that $\cH_i
:= \{H_{i,1}, H_{i,2}, \ldots \}$ is a generating set for $\cP_i$. But
$\cH_i \subseteq \cQ_i$, so $\cP_i = [\cH_i] \subseteq \cQ_i$. 

By the same reasoning, there is a permutation $\tau$ such that $\cQ_i
\subseteq \cP_{\tau(i)}$. We cannot relabel the $\cP_i$'s as well, but
if $\tau^k(i) = i$, then we have $\cP_i \subseteq \cQ_i \subseteq
\cP_{\tau(i)} \subseteq \cQ_{\tau(i)} \subseteq \cP_{\tau^2(i)}
\subseteq \cQ_{\tau^2(i)} \subseteq \cdots \subseteq \cP_{\tau^k(i)} =
\cP_i$, so we must have equality throughout; in particular, $\cP_i =
\cQ_i$ for each $i$. $\hfill\square$

\medskip

Now for the proof of Theorem \ref{m=n}.

%\resnum{Theorem\label{m=n}}
%{Let 
%$\cP_1 \circ \cdots \circ \cP_m$ be a factorisation of the
%hereditary compositive property $\cP$ into indecomposable
%hereditary compositive properties. Then $m = dc(\cP)$. 
%}

\medskip\noindent{\bf Proof of Theorem \ref{m=n}:}
Given any generating set $\cG_i$ for $\cP_i$, every graph in $\cG_1 +
\cdots + \cG_m$ has decomposability $m$ by construction. Then every
graph in $(\cG_1 + \cdots + \cG_m)^{\td} \subseteq \cM^*(\cP)$ has
decomposability at least $m$, so $dc(\cP) \geq m$. 

%To prove the reverse inequality, we suppose $m < n := dc(\cP)$ and then construct
%generating sets  for $\cP$, until we get a contradiction. 
If $m < n := dc(\cP)$, and $G$ is a $\cP$-maximal graph with decomposability $n$, then in any $(\cP_1, \ldots, \cP_m)$-partition of $G$ some $\cP_i$-part is the join of  two or more ind-parts. There is only a finite number of ways in which this can happen, and we will construct a sequence of generating sets so that each one excludes at least one of the possibilities until we reach a contradiction.

When graphs or sets have a double subscript, we will
use the second number to denote which step of our construction we are
in. For each $i$, we start with some generating set $\cG_i$ consisting
only of indecomposable $\cP_i$-strict graphs. 
\newline

Let $H_1 \in (\cG_1 + \cdots + \cG_m)^{\td}$; then $H_1$ is a join
$H_{1,1} + \cdots + H_{n,1}$ of $n$ ind-parts. In general suppose we
have graphs $H_1, H_2, \ldots, H_{k-1}$ such that, for each $s=1,2,\dots,k-1$: 
\begin{itemize}
\item[(a)] $H_s$ is $\cP$-maximal;
\item[(b)] $dc(H_s) = n$, and  $H_s = H_{1,s} + \cdots +
H_{n,s}$; 
\item[(c)] for $j = 1,\ldots,n$, $H_{j,1} \subseteq H_{j,2} \subseteq
\cdots \subseteq H_{j,k-1}$; and
\item[(d)] there is a partition $(J_{1,s},J_{2,s},\dots,
J_{m,s})$ of $\{1,2,\ldots,n\}$ such that $$ \sum_{j \in
J_{i,s}} H_{j,s} \in \cP_i.$$ 
\end{itemize}

Now let $H_{k} \in (\cG_1[\sum_{j \in J_{1,(k-1)}} H_{j,(k-1)}] + \cdots
+ \cG_m[\sum_{j \in J_{m,(k-1)}} H_{j,(k-1)}])^{\td}$.  As 
$H_k$ contains $H_{k-1}$, by Lemma~\ref{gen-0} we can label the
ind-parts of $H_k = H_{1,k} + \cdots + H_{n,k}$ so that $H_{1,(k-1)}
\subseteq H_{1,k}, \ldots, H_{n,(k-1)} \subseteq H_{n,k}$. It is
important to note that the indecomposable graph from $\cG_i[\sum_{j
\in J_{i,(k-1)}} H_{j,(k-1)}]$ therefore spans $\sum_{j \in 
J_{i,(k-1)}} H_{j,k}$ (note the change in subscript).  
By Lemma~\ref{max-char} there is a partition $(J_{1,k}, \cdots,
 J_{m,k})$ of $\{1,2,\ldots,n\}$ so that $\sum_{j \in J_{i,k}}
H_{j,k} \in P_i$. 
\newline

Since there is only a finite number of partitions of
$\{1,2,\ldots,n\}$, at some step $B$ we must end up with a partition
that occurred at some previous step $A < B$. Without loss of
generality, suppose that $|J_{1,A}| = r \geq 2$. Then $\sum_{j\in
J_{1,A}}H_{j,A}  \in \cP_1$; the
indecomposable graph $H_A$ from $\cG_1[\sum_{j\in J_{1,A}}
H_{j,A}]$ that is used in step $A+1$ spans $\sum_{j\in
J_{1,A}}H_{j,(A+1)} 
$; this join properly contains the $\cP_1$-maximal graph
$H_A$ and  therefore is not in $\cP_1$. But, for each $j$, $H_{j,(A+1)}
\subseteq H_{j,(A+2)} \subseteq \cdots \subseteq H_{j,B}$, and so
$\sum_{j\in J_{1,A}}H_{j,(A+1)} \subseteq \sum_{j\in J_{1,A}}H_{j,B}$.
But $J_{1,A}=J_{1,B}$ and $\sum_{j\in J_{1,B}}H_{j,B}
\in \cP_1$, so $\sum_{j\in J_{1,A}}H_{j,(A+1)}\in\cP_1$,  a contradiction. 

Thus we must have $|J_{i,A}| = 1$, for each $i = 1,2,\ldots,m$, and so $m=n$.
$\hfill\square$

\section{Definitions and results from~\cite{uft-2}\label{sec-ind-def}}

This section and the next are the induced-hereditary analogues of
Sections~\ref{sec-def} and~\ref{sec-proof}, along with a highly
important result (Theorem~\ref{m-resp-n}) adapted from~\cite{uft-2}.  

In~\cite{uft-2} Mih\'{o}k generalised the results of~\cite{uft-1} from
additive hereditary properties to the wider class of 
additive~\emph{induced}-hereditary graph properties (we point out again, 
though, that hereditary \emph{compositive} properties
are not all additive); the concepts introduced in that article are presented
here. We caution the reader that there are some significant
differences between the old and new definitions of ``generating set'',
``join'', ``decomposability'', ``$\cP$-strict'' and ``ind-part'';
these new definitions will apply throughout the rest of the paper,
even for hereditary properties (that are \emph{a fortiori}
induced-hereditary). 

We use $G \leq H$ to denote that $G$ is an induced subgraph of $H$. 
The smallest induced-hereditary property that contains a set $\cG$ is 
denoted by $\la \cG\ra$.  This is the  induced-hereditary property 
{\em generated by $\cG$\/}, or that {\em $\cG$ generates\/}. We say that 
$\cG$ {\em is a generating set for\/}  $\cP$ if $\la \cG\ra = \cP$. 
It is easy to see that:

\[\begin{array}{lll}
\la \cG\ra &=& \{G \mid \exists\, H \in \cG,\ G \leq H \}.
\end{array}\]

The {\em $*$-join\/} of $n$ graphs $G_1, \ldots, G_n$ is the set
\[ G_1 * \cdots * G_n := \{G \mid \bigcup_{i=1}^n G_i \subseteq G
\subseteq \sum_{i=1}^n G_i \} \] 
where $\bigcup$ and $\sum$ represent disjoint union and join, respectively.
Given $n$ sets of graphs, we define their $*$-join by 
\[ S_1 * \cdots * S_n := \bigcup \left(G_1 * \cdots * G_n\right)\ ,\] the
union being over all ways of the selecting the $G_i$ so that $G_i \in
S_i$ for all $i$.
We note that this is just the same as $S_1 \circ \cdots \circ S_n$,
but it is aesthetically pleasing to have the $*$ notation. If $\cP_1,
\ldots, \cP_n$ are additive properties, and $G_i \in \cP_i$ for all
$i$, then for all positive integers $k$ we have 
\[ kG_1 * \cdots * kG_n \subseteq \cP_1 \circ \cdots \circ \cP_n \] 
where $kG$ is the disjoint union of $k$ copies of $G$.

A {\em $\cP$-decomposition of $G$ with $n$ parts\/} is a partition
$(V_1,\ldots,V_n)$ of $V(G)$ such that for all $i$ $V_i \ne
\emptyset$, and for all positive integers $k$ we have $kG[V_1] *
\cdots * kG[V_n] \subseteq \cP$. The {\em $\cP$-decomposability
number\/} $dec_{\cP}(G)$ of $G$ is the maximum number of parts in a
$\cP$-decomposition of $G$; if $G \not\in \cP$, then we put $dec_{\cP}(G) =
0$. If $G\in \cP$, then, for all positive integers $k$, $kG\in\cP$;
therefore $G\in\cP$ if and only if $dec_{\cP}(G)\ge 1$.  
Also, $G$ is {\em $\cP$-decomposable\/} if $dec_{\cP}(G) >1$. 
If $\cP$ is the product of two additive induced-hereditary properties, then
\emph{every} graph in $\cP$ with at least two vertices is
$\cP$-decomposable. 

\stmtnum{Lemma\label{ind-her-0}}
{Let $\cP = \cP_1 \circ \cdots \circ \cP_m$, where the $\cP_i$'s are
additive properties of graphs. Then any $(\cP_1, \ldots,
\cP_m)$-partition of a graph $G$ is a $\cP$-decomposition of $G$. If
the $\cP_i$'s are induced-hereditary, then every graph in $\cP$ with at least
$m$ vertices has a partition with all $m$ parts non-empty. 
}
\stmtnumend
\newline

A graph $G$ is {\em $\cP$-strict\/} if $G \in \cP$ but $G * K_1 
\nsubseteq \cP$; we denote the set of $\cP$-strict graphs by
$\bS(\cP)$.  
If $f(\cP) = \min\{|V(F)| \mid F \not\in \cP \}$, then $G * K_1 * \cdots
* K_1 \nsubseteq \cP$, where the $*$ operation is repeated $f(\cP)$
times. Thus, every $G \in \cP$ is an induced-subgraph of some
$\cP$-strict graph (with fewer than $|V(G)| + f(\cP)$ vertices), and so
$\la \bS(\cP)\ra = \cP$. 
Similarly, $dec_{\cP}(G) < f(\cP)$.  

The {\em decomposability number\/} $dec(\cG)$ 
of a generating set $\cG$ of $\cP$
is $$\min\{dec_{\cP}(G) \mid  G \in \cG \};$$ the
{\em decomposability number\/} $dec(\cP)$ of $\cP$ is 
$dec(\bS(\cP))$. A property with $dec(\cP) = 1$ is 
{\em indecomposable\/}. An indecomposable property is also irreducible
and it will turn out that the converse is also true. 

\stmtnum{Lemma\label{ind-her-1}}
{Let $\cP_1, \ldots, \cP_m$ be induced-hereditary properties of
graphs, and let $G$ be a  $\cP_1 \circ \cdots \circ \cP_m$-strict
graph. Then, for every $(\cP_1, \ldots, \cP_m)$-partition $(V_1,
\ldots, V_m)$ of $V(G)$, $G[V_i]$ is $\cP_i$-strict (and in particular
non-empty). 
}
%If $G[V_1]*K_1 \subseteq \cP_1$, then $G * K_1 \subseteq 
%(G[V_1]*K_1) * G[V_2] * \cdots * G[V_m] \subseteq 
%\cP_1 \circ \cdots \circ \cP_m$. 
\stmtnumend
\newline

It follows that $dec(\cA \circ \cB) \geq dec(\cA) + dec(\cB)$, and
thus any factorisation of an additive induced-hereditary property $\cP$ has 
at most $dec(\cP)$ irreducible additive induced-hereditary factors.

\resnum[\ \cite{uft-2}]{Lemma\label{ind-her-2}}
{Let $G$ be a $\cP$-strict graph, and let  $G' \in \cP$ be an induced
supergraph of $G$, i.e., $G \leq G'$. Then $G'$ is $\cP$-strict, and
$dec_{\cP}(G) \geq dec_{\cP}(G')$. 
}
Every graph in $G * K_1$ is an induced subgraph of a graph in $G' *
K_1$, so $G'$ must be $\cP$-strict. If $(V_1, \ldots, V_n)$ is a
$\cP$-decomposition of $G'$ with $n$ parts, then $(V_1 \cap V(G),
\ldots, V_n \cap V(G))$ is a $\cP$-decomposition of $G$; moreover, it
has $n$ parts unless, for some $i$, $V_i \cap V(G) = \emptyset$, which
is impossible because $G$ is $\cP$-strict. 
\resnumend
%\newline

\stmtnum[\ \cite{uft-2}]{Lemma\label{ind-her-3}}
{If $\cG$ generates the induced-hereditary property $\cP$, 
then $dec(\cG) \leq dec(\bS(\cP))$, with equality if $\cG \subseteq
\bS(\cP)$.  
}
\stmtnumend
\newline

For $\cG \subseteq \cP$, and $H \in \cP$, let $\cG[H] := \{G \in \cG \mid
H \leq G \}$.  

\stmtnum[\ \cite{uft-2}]{Lemma\label{ind-her-4}}
{Let $\cG$ generate the additive induced-hereditary property $\cP$,
and let $H$ be an arbitrary graph in $\cP$. Then $\cG[H]$ also
generates $\cP$. 
}\stmtnumend
%For any graph $G' \in \cP$, the graph $G' \cup H$ is also in $\cP$ (by
%additivity). Since $\cG$ is a generating set, there is some graph $H'
%\in \cG$ that contains $G' \cup H$ as an induced-subgraph; but then
%$G' \leq H' \in \cG[H]$. 
%\resnumend
\newline

For a generating set $\cG$, let $\cG^{\td} := \{G \in \cG \mid G \in
\bS(\cP),\ dec_{\cP}(G) = dec(\bS(\cP)) \}$. The following is a simple
consequence of Lemmas \ref{ind-her-2} and \ref{ind-her-4}.

\stmtnum[\ \cite{uft-2}]{Lemma\label{ind-her-5}}
{If $\cG$ generates the additive induced-hereditary property $\cP$,
then so does $\cG^{\td}$. 
}\stmtnumend
%Let $H \in \bS(\cP)$ be a graph with $dec_{\cP}(H) =
%dec(\bS(\cP))$. By Lemma~\ref{ind-her-2} $\cG[H] \subseteq \cG^{\td}$,
%and by 
%Lemma~\ref{ind-her-4} $\cG[H]$ is a generating set for $\cP$.
%\resnumend
\newline

A graph $G\in\cP$ is {\em uniquely $\cP$-decomposable\/} if there is only one
$\cP$-dec\-omp\-osition of $G$ with $dec_{\cP}(G)$ parts. If $\cP =
\cP_1 \circ \cdots \circ \cP_n$, then by Lemma~\ref{ind-her-0} a
uniquely $\cP$-decomposable graph $G$ with $dec_{\cP}(G) = n$ must be
uniquely $\{\cP_1, \ldots, \cP_n\}$-partition\-able (every $\{\cP_1,
\ldots, \cP_n\}$-partition gives the same unordered partition of
$V(G)$). 
If $(V_1, \ldots, V_n)$ is the unique $\cP$-decomposition of $G$, we
call the graphs $G[V_1], \ldots, G[V_n]$ its {\em ind-parts} (although they
are themselves usually $\cP$-dec\-omposable).

\stmtnum{Lemma\label{ind-her-6}}
{Let $G$ be a graph in $\bS(\cP)$ with $dec_{\cP}(G) = dec(\cP)$, and
suppose that $G$ has a unique $\cP$-decomposition $(V_1, \ldots,
V_{dec(\cP)})$ with $dec(\cP)$ parts. If $G \leq H$, then $H \in
\bS(\cP)$, $dec_{\cP}(H) = dec(\cP)$, and, for any $\cP$-decomposition
$(W_1, \ldots, W_{dec(\cP)})$ of $H$, we can relabel the $W_i$'s so
that, for all $i$, $W_i \cap V(G) = V_i$. 
}
\stmtnumend
\newline

In the hereditary case it was very important that if $G = G_1 + \cdots
+ G_m$, each $G_i$ is the join of ind-parts (the partition into
ind-parts ``respected'' the partition into $G_i$'s); in the
induced-hereditary case we can prove analogous results
(e.g., Corollaries~\ref{respect}--\ref{unique-respect}) for
$\cP$-strict, uniquely $\cP$-decomposable graphs with $dec_{\cP}(G) =
dec(\cP)$, which allows us to generalise Theorem~\ref{m=n}.  

\defnnum{Definitions}
{Let $d_0 = (U_1, U_2, \ldots, U_m)$
 be a $\cP$-decomposition of a graph $G$.
 A $\cP$-decomposition $d_1=(V_1,V_2,\dots,V_n)$ of $G$ {\em respects\/} $d_0$ 
if no $V_i$ intersects two or more
$U_j$'s; that is, each $V_i$ is contained in some $U_j$, and so each
$U_j$ is a union of $V_i$'s. 

If $G$ is a graph, then $s\cstar G$ denotes the set $G*G*\cdots*G$, where
there are $s$ copies of $G$.  For $G^* \in s\cstar G$,
denote the copies 
of $G$ by $G^1, \ldots, G^s$. Then {\em $G^*$ respects $d_0$\/} if
$G^* \in s G[U_1] * \cdots * s G[U_m]$; that is,  
two vertices in different $G^i$'s are joined by an edge only if they
are also contained in different $U_j$'s. A $\cP$-decomposition 
$d = (V_1, \ldots, V_n)$ of $G^*$ {\em respects $d_0$
uniformly\/} if, for each $V_i$, there is a $U_j$ such that, for every $G^k$, 
$V_i \cap V(G^k)\subseteq U_j$.  
The decomposition of $G^k$ induced by $d$
is denoted $d|G^k$.

If $G$ is uniquely $\cP$-decomposable, its ind-parts {\em respect
$d_0$\/} if its unique $\cP$-decomposition with $dec_{\cP}(G)$ parts
respects $d_0$. If $G^*$ is uniquely $\cP$-decompos\-able, its
ind-parts {\em respect $d_0$ uniformly\/} if: (a) for some $s$,
 $G^* \in s\cstar
G$; (b) $G^*$ respects $d_0$; and (c) $G^*$'s unique
$\cP$-decomposition with $dec_{\cP}(G^*)$ parts respects $d_0$
uniformly. 

\bigskip

\begin{figure}[htb]
\begin{center}
\input{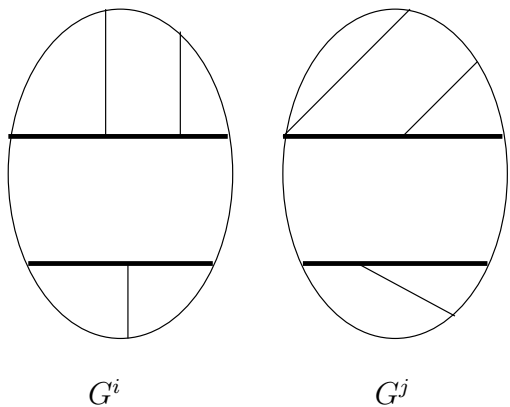} 
\caption{$d$ (vertical lines) respects $d_0$ (horizontal lines) uniformly}
\label{Fig-uniform-respect}
\end{center}
\end{figure}

\begin{figure}[htb]
\begin{center}
\input{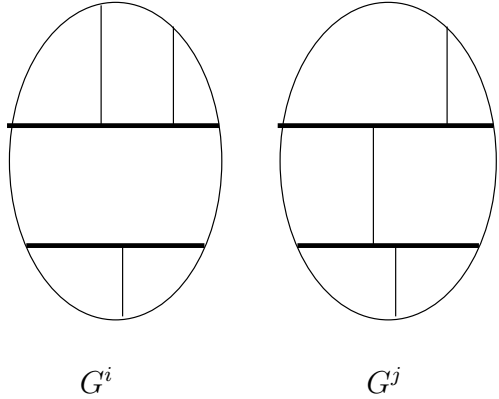} 
\caption{$d$ (vertical lines) respects $d_0$ (horizontal lines) on
both $G^i$ and $G^j$, but not uniformly} 
\label{Fig-respect}
\end{center}
\end{figure}

\bigskip

The {\em extension of $d_0$\/} to $G^*$ is the decomposition obtained
by repeating $d_0$ on each copy of $G$. If $G^*$ respects $d_0$, or if
it has a $\cP$-decomposition that respects $d_0$ uniformly, then the
extension of $d_0$ is also a $\cP$-decomposition of $G^*$. 
In particular, $G^*$ is a graph in $\cP$.

We will sometimes write $G^i \cap U_x$ (or just $U_x$ when it is clear
we are referring to $G^i$) to mean the vertices of $G^i$ that
correspond to $U_x$, and $G^* \cap U_x$ (or just $U_x$, when it is
clear from the context) to mean $\bigcup_i (G^i \cap U_x)$. 
}
\defnnumend
\newline

The required result is a corollary of the following theorem of
Mih\'{o}k; he actually proved it when $m = n$
(Corollary~\ref{unique-super}), but very little modification is needed
to establish the general case, and we follow his proof and notation
rather closely. 

\resnum{Theorem\label{m-resp-n}}
{Let $G$ be a $\cP$-strict graph with $dec_{\cP}(G) = n$, and let $d_0
= (U_1, U_2, \ldots, U_m)$ be a fixed $\cP$-decomposition of
$G$. Then there is a $\cP$-strict graph $G^* \in s\cstar G$ (for some
$s$) that respects $d_0$, and moreover any $\cP$-decomposition of
$G^*$ with $n$ parts respects $d_0$ uniformly. 
}
Let $d_i = (V_{i,1},V_{i,2}, \ldots ,V_{i,n}),\  i=1, \ldots, r$, be
the  $\cP$-decompositions of $G$ with $n$ parts which do not respect
$d_0$. Since $G$ is a finite graph, $r$ is a nonnegative integer. If
$r=0$, take $G^* = G$; otherwise we will construct a  graph $G^* =
G^*(r)\in s\cstar G$ as above, denoting the $s$ copies of $G$ by
$G^1, \ldots, G^s$. 

If the resulting $G^*$ has a $\cP$-decomposition $d$ with $n$ parts,
then, since $G$ is $\cP$-strict, $d|G^i$ will also have $n$ parts. The
aim of the construction is to add new edges $E^* = E^*(r)$ to
$sG$ to exclude the possibility that $d|G^i = d_j$, for any $1 \leq
i \leq s, 1 \leq j \leq r$. We will only add edges between $G^i
\cap U_x$ and $G^j \cap U_y$, where $i \not= j$ and $x \not= y$, so
that $G^*$ will respect $d_0$. 

% Explain that G^* is thus in P.

We shall use two types of constructions.

\vspace{5mm}
{\noindent\bf Construction 1.} $G^i \Rightarrow G^j$.  

This is a
graph in $2\cstar G$ such that, if $d$ is a $\cP$-decomposition of
$G^i\Rightarrow G^j$ and $d|G^i$ respects $d_0$, then $d|G^j$
respects $d_0$; moreover, $d$ respects $d_0$ uniformly on
$G^i\Rightarrow G^j$.  (We comment that this corrects a minor error in
\cite{uft-2}.  The author of \cite{uft-2} was independently
 aware of both the error and its correction.)

Since $G$ is $\cP$-strict, 
there is a graph  $F \in (G * K_1)\setminus \cP$. Let $N_F(z)$ 
be the neighbours in $G$ of $z \in V(K_1)$. For $y=1,2,\dots,m$, 
let $Z_y$ denote $U_y \cap
N_F(z)$.
Let $G^i,G^j, i \not = j$ be disjoint copies of $G$; join every vertex
of $U_x$ in  $G^j$ to every vertex of $Z_y$, $x \not = y,$
in  $G^i$.  Note that $G^i\Rightarrow G^j\in 2G[U_1]*2G[U_2]*\cdots
*2G[U_m]$.  Since $d_0$ is a $\cP$-decomposition of $G$, this implies
that $(G^i\Rightarrow G^j)\in \cP$.

Let $d=(V_1,V_2, \dots, V_\ell)$ be a $\cP$-decomposition of
$H=(G^i\Rightarrow G^j)$ such that $d|G^i$
respects $d_0$, but $d|G^j$ does not respect $d_0$ 
(or at least, not in the same manner, i.e., 
$d$ does not respect $d_0$ uniformly).
Then there is a $k$ such that $V_k \cap G^i \subseteq U_y$, but $v \in
V_k \cap G^j$ belongs to $U_x$, $x \ne y$. 
We claim $F$ is an induced subgraph of a graph in 
$ H[V_1]*H[V_2]*\cdots*H[V_\ell]$, which implies $F\in\cP$, a contradiction.

To see this, consider the vertex $v$ and $G^i$.  We have edges from
$v$ to every vertex in $Z_w$, for all $w\ne x$.  We are only missing
the edges from $v$ to $Z_x\cap G^i$.  But $d|G^i$ respects $d_0$, so
$U_x\cap G^i$ is the union of, say, $V_{t_1}\cap G^i, V_{t_2}\cap
G^i,\dots, V_{t_r}\cap G^i$. Since $V_k\cap G^i\subseteq U_y$ and
$y\ne x$, $k\notin \{t_1,t_2,\dots,t_r\}$.  Since $d$ is a
$\cP$-decomposition of $H$, we may freely add edges between 
$V_i$'s and remain in $\cP$; in
particular, one graph in $H[V_1]*H[V_2]*\cdots*H[V_\ell]$ is obtained
by adding  precisely the edges between $v$ and $Z_x\cap
G^i$.  Clearly  $F$ is the subgraph of this graph induced by $G^i$
and $v$, as claimed.

%$\cP$-decomposition of $H = G^i \Leftrightarrow G^j$ since there is a
%graph in  
%$H[V_1] * H[V_2] * \cdots * H[V_n]$ 
%which contains an induced copy of $F$ (we can add the appropriate
%edges between $v$ and $Z_x$). 

\vspace{5mm}
{\noindent\bf Construction 2.} $m\bullet k_t G$.  

For a $\cP$-decomposition 
$d_t=(V_{t, 1},V_{t, 2},\dots,V_{t, dec_{\cP}(G)})$
 of $G$ that does not respect $d_0$,
$m\bullet k_t G$ is a
graph in $(mk_t)\cstar G$ having no $\cP$-decomposition
$d=(W_1,W_2,\dots,W_{dec_{\cP}(G)})$ such that, for all of the $mk_t$
induced copies $G^i$ of $G$, $d|G^i = d_t$.

\begin{figure}[htb]
%+%\hspace{-10mm}\input{Constr2.pstex_t} 
\hspace{-10mm}\input{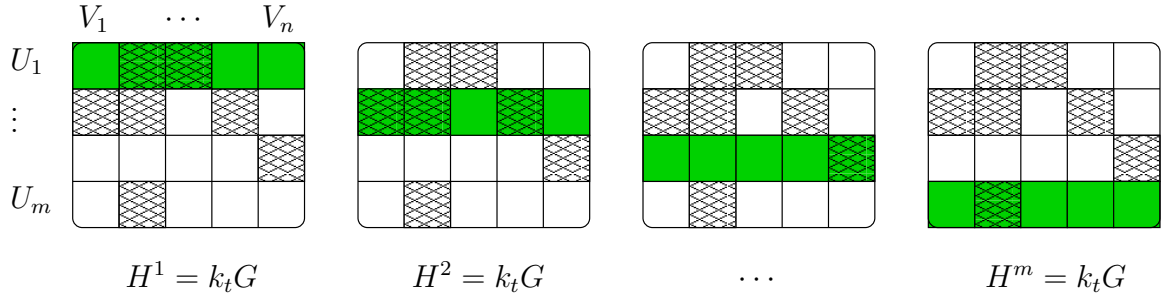} 
\caption{$m \bullet k_t G$ - we only put edges between the $m$ shaded parts}
\label{Fig-constr-2}
\end{figure}

Let $n=dec_{\cP}(G)$ and
let  $A_{i,j}(t)$ denote $U_i \cap V_{t,j}, 1 \leq i
\leq m, 1 \leq j \leq n$. Since $d_t$ does not respect $d_0$, at
least $n+1$ sets $A_{i,j}(t)$ are nonempty. Because $dec_{\cP}(G)=n$,
there exists a positive integer $k_t$ such that $k_t
G[A_{1,1}(t)] * k_t G[A_{1,2}(t)] * \cdots * k_t
G[A_{m,n}(t)] \not \subset \cP$. Fix a graph  $F_t \in (k_t
G[A_{1,1}(t)] * k_t G[A_{1,2}(t)] * \cdots * k_t
G[A_{m,n}(t)])\setminus \cP$. 
%Denote by $E_{ij;i'j'}(t)$ the
%set of edges of $F_t$ joining the vertices of  $k_t
%G[A_{i,j}(t)]$ and $k_t G[A_{i',j'}(t)]$. 

%Before constructing $m.k_t \bullet G$ we explain the idea behind the 
%construction. 
Suppose that in $H = k_t G$ we replace the edges between $H \cap U_x$ and $H \cap U_y$ by the edges that there are between $F_t \cap U_x$ and $F_t \cap U_y$, for all $x \not= y$; the $U_x$'s still form a $\cP$-decomposition of the resulting graph, $\tilde{H}$, so it is in $\cP$.
If the extension of $d_t$ were also a $\cP$-decomposition of $\tilde{H}$, we could obtain $F_t$ immediately by replacing the edges between $\tilde{H} \cap V_{t, i}$ and $\tilde{H} \cap V_{t,j}$, by those between $F_t \cap V_{t, i}$ and $F_t \cap V_{t,j}$, for all $i \not= j$. The only problem is that $\tilde{H}$ does not contain $k_t$ disjoint copies of $G$, as we altered edges \emph{inside} the copies of $G$.

So instead we construct $m \bullet k_t G$ from $m$ disjoint
copies of $H = k_t G$, denoted by $H^j, j=1,2, \ldots, m$ (see 
Figure~\ref{Fig-constr-2}). We add 
edges between $H^x \cap U_x$ and $H^y \cap U_y$, corresponding to 
the edges that there are between $F_t \cap U_x$ and $F_t \cap U_y$, 
for all $x \not= y$.
Because $d_0=(U_1,U_2,\dots,U_m)$ is a $\cP$-decomposition, $H$ is in $\cP$.

Now $H^1 \cap U_1, \ldots, H^m \cap U_m$ form a copy of $\tilde{H}$ in
$m \bullet k_t G$. Suppose 
$H' = m\bullet k_t G$ has a $\cP$-decomposition
$d=(W_1,W_2,\dots,W_{n})$ such that, for every one of the
$mk_t$ 
induced copies $G^i$ of $G$, $d|G^i=d_t$; then we can obtain $F_t$ 
as an induced subgraph of a graph in $H'[W_1]*H'[W_2]*\cdots*H'[W_n]$
(by changing edges in the copy of $\tilde{H}$ as explained above). 

\begin{figure}[htb]
\begin{center}
\input{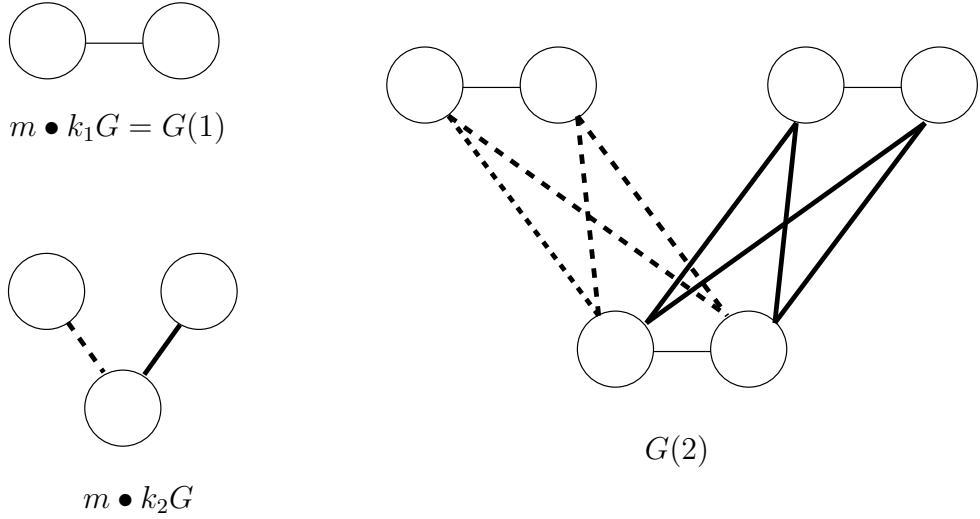} 
\caption{Constructing $G(2)$ from $G(1)$ and $m \bullet k_2 G$}
\label{Fig-constr-G-2}
\end{center}
\end{figure}

\vspace{5mm}
We now construct $G^*$ as follows.
First let $G(1) := m \bullet k_1 G$. For $1 < \ell \leq r$, construct
$G(\ell)$ by taking  $m k_\ell$ disjoint copies
$G(\ell-1)^1,\dots,G(\ell-1)^{mk_\ell}$ of $G(\ell-1)$.  For each copy $G^i$ of
$G$ in $G(\ell-1)^i$ and each copy $G^j$ of $G$ in $G(\ell-1)^j$, we add the
edges between them that are between the $i^{{th}}$ and
$j^{{th}}$ copies of $G$ in $m\bullet k_\ell G$.  
(See Figure~\ref{Fig-constr-G-2}.)

\bigskip

\bigskip

Finally, from $G(r)$, which is in, say, $s\cstar  G$, consisting of copies
$G^1$, $G^2$, \dots, $G^s$ of $G$, we create $G^*$ by adding 
two more copies $G^0$ and
$G^{s+1}$ of $G$.  For each $i\in\{1,2,\dots,s\}$, we add the edges
between $G^0$ and $G^i$ to create the graph $G^0\Rightarrow G^i$, we
add the edges between $G^i$ and $G^{s+1}$ to create the graph
$G^i\Rightarrow G^{s+1}$, and we add the edges between $G^{s+1}$ and
$G^0$ to create the graph $G^{s+1}\Rightarrow G^0$.

Let $d$ be a $\cP$-decomposition of $G^*$ with $n$ parts (it might be
that none exists, in which case we are done).  
For $1 \leq \ell \leq r$, if every copy of $G(\ell-1)$ in $G(\ell)$ contains
a copy of $G$ for which $d|G = d_\ell$, then we would have $m k_\ell$ such
copies of $G$ inducing a copy of $m \bullet k_\ell G$, which we know is
impossible. So by induction from $r$ to $1$, there is a copy $G^p$
of $G$ for which $d|G^p$ is none of $d_1,d_2,\dots,d_r$.  Thus,
$d|G^p$ respects $d_0$.  But $G^p\Rightarrow G^{s+1}$ is an induced
subgraph of $G^*$, so $d|G^{s+1}=d_0$ (and in fact $d$ respects $d_0$
uniformly on these two copies of $G$).  Similarly, $d|G^0$ respects
$d_0$ and, again in the same way, $d$ respects $d_0$ uniformly, as
required.  
\resnumend
%\newline

\resnum{Corollary\label{respect}}
{Let $G$ be a $\cP$-strict graph with $dec_{\cP}(G) = dec(\cP)$, and
let $d_0 = (U_1, U_2, \ldots, U_m)$ be a fixed $\cP$-decomposition
of $G$. Then there is a $\cP$-decomposition of $G$ with exactly
$dec(\cP)$ parts that respects $d_0$. 
}
In Theorem~\ref{m-resp-n}, since $G^* \geq G$ we know $G^*$ is
$\cP$-strict, and so $dec(\cP) \leq dec_{\cP}(G^*) \leq dec_{\cP}(G) =
dec(\cP)$. Thus $G^*$ has at least one $\cP$-decomposition $d$ with
$dec(\cP)$ parts; $d|G$ also has $dec(\cP)$ parts (since $G$ is
$\cP$-strict) and respects $d_0$. 
\resnumend
%\newline

\resnum[\ \cite{uft-2}]{Corollary\label{unique-super}}
{Let $G$ be a $\cP$-strict graph with $dec_{\cP}(G) = n$, and let $d_0
= (U_1, U_2, \ldots, U_n)$ be a fixed $\cP$-decomposition of $G$ with
$n$ parts. Then there is a $\cP$-strict graph $G^* \in s\cstar G$ (for
some $s$) which has a unique $\cP$-decomposition $d$ with $n$ parts,
and $d|G^j = d_0$ for all $j$. 
}
The only $\cP$-decomposition of $G$ with $n$ parts that respects $d_0$
is $d_0$ itself (since here $d_0$ has exactly $n$ parts). Thus in
Theorem~\ref{m-resp-n}, the only possible decomposition of $G^*$ with
$n$ parts is the extension of $d_0$, which is a $\cP$-decomposition of
$G^*$ by construction. 
\resnumend 
%\newline

\medskip

The next result tells us that under certain conditions, given a factorisation $\cQ_1 \circ \cdots \circ \cQ_m$ of $\cP$ into additive induced-hereditary properties, and a $\cP$-decomposition $d_0$ of $G$, we can group the parts of $d_0$ to get a $(\cQ_1, \ldots, \cQ_m)$-partition of $G$. Of course, $d_0$ does not respect \emph{all} $(\cQ_1, \ldots, \cQ_m)$-partitions; in fact, if $m = dec(\cP)$, $d_0$ can only respect one partition, namely $d_0$ itself (note that none of the parts of a partition can be empty, because $G$ is $\cP$-strict). 
We will see later (Theorem~\ref{ind-her-m=n}) that when we factor the $\cQ_i$'s as far as possible we do get exactly $dec(\cP)$ irreducible factors, say $\cP_1, \ldots, \cP_{dec(\cP)}$, and applying the corollary we get that $d_0$ \emph{is} a $(\cP_1, \ldots, \cP_{dec(\cP)})$-partition.
%We will actually use Corollary~\ref{colour-respect-1} to show that when we factor 
%the $\cQ_i$'s as far as possible we do get exactly $dec(\cP)$ irreducible factors 
%(Theorem~\ref{ind-her-m=n}), say $\cP_1, \ldots, \cP_{dec(\cP)}$, and applying the 
%corollary we get that $d_0$ \emph{is} a $(\cP_1, \ldots, \cP_{dec(\cP)})$-partition.

\resnum{Corollary\label{colour-respect-1}}
{Let $G$ be a $\cP$-strict graph with $dec_{\cP}(G) = dec(\cP)$, and with some $\cP$-decomposition $d_0 = (U_1, U_2, \ldots, U_{dec(\cP)})$. If $\cP = \cQ_1 \circ \cdots \circ \cQ_m$, then $G$ has a $(\cQ_1, \ldots, \cQ_m)$-partition that $d_0$ respects.
}
The graph $G^*$ of Corollary~\ref{unique-super} has 
some $(\cQ_1, \ldots, \cQ_m)$-partition $d_1$; this is also a $\cP$-decomposition. By Corollary~\ref{respect} the unique $\cP$-decomposition $d$ of $G^*$ with $dec(\cP)$ parts must respect $d_1$; and the restriction of $d$ to $G$ is just $d_0$.
\resnumend
\newline

The set of $\cP$-strict,
uniquely $\cP$-decomposable graphs with $dec_{\cP}(G) = dec(\cP)$ is
denoted $\bS^{\tD}(\cP)$, or just $\bS^{\tD}$. By
Lemma~\ref{ind-her-5} and Corollary~\ref{unique-super} 
$\bS^{\tD}$ is a generating set for $\cP$; in fact, for any $G \in
\bS^{\td}$ and any specific $\cP$-decomposition $d$ of $G$, we can
find an induced supergraph in $\bS^{\tD}$ whose ind-parts uniformly
respect $d$. 

\stmtnum{Corollary\label{unique-respect}}
{Let $G$ be a $\cP$-strict graph with $dec_{\cP}(G) = dec(\cP)$, and
let $d_0 = (U_1, U_2, \ldots, U_m)$ be a fixed $\cP$-decomposition
of $G$. Then there is a uniquely $\cP$-decomposable $\cP$-strict graph
$G^* \geq G$ whose ind-parts respect $d_0$ uniformly. 
}
\stmtnumend
%\newline

\resnum{Corollary\label{dec-part}}
{Let $\cP = \cP_1 \circ \cdots \circ \cP_{dec(\cP)}$.
Let $G$ be a $\cP$-strict graph with
$dec_{\cP}(G) = dec(\cP)$.  If 
$d_0=(U_1,U_2,\dots,U_m)$ is a $\cP$-decomposition of $G$, then there
is a factorization $\cP =\cQ_1 \circ \cdots \circ \cQ_m$ such that 
$d_0$ is a $(\cQ_1, \ldots, \cQ_m)$-partition of $G$. 
}
By Corollary~\ref{unique-respect} there is a a uniquely
$\cP$-decomposable graph $G^* \geq G$ whose ind-parts respect
$d_0$ uniformly.  Let $(V_1,V_2,\dots,V_{dec(\cP)})$ be the unique
$\cP$-decomposition of $G^*$.
By Lemma~\ref{ind-her-0}, 
the ind-parts of $G^*$ must form its unique $(\cP_1, \ldots,
\cP_{dec(\cP)})$-partition, so there is a partition $(J_1,J_2,\dots,J_m)$ of
$\{1,2,\dots,dec(\cP)\}$ such that, for each $i$, $\displaystyle U_i =
\cup_{j\in J_i}V_j$ (when we restrict the $V_j$ to a particular copy of $G$
in $G^*$).  It follows that $\displaystyle G[U_i]\in \prod_{j\in
J_i}\cP_j$, so we may set $\displaystyle \cQ_i=\prod_{j\in
J_i}\cP_j$.  
\resnumend
%\newline

\section{Unique factorisation for additive induced-hereditary
properties\label{sec-ind-proof}} 
The strategy for proving the uniqueness of the factorisation of an
additive induced-hereditary property into irreducible additive
induced-hereditary properties is the same as for hereditary compositive
properties.  We shall first show that there is at most one into
$dec(\cP)$ factors and then that any such factorisation must have
$dec(\cP)$ factors.

The following construction of a generating set for $\cP$ will be
essential in proving unique factorisation. Suppose we are given a
factorisation $\cP = \cP_1 \circ \cdots \circ \cP_m$ into
indecomposable additive induced-hereditary factors, and, for each $i$,
we are given a generating set $\cG_i$ of $\cP_i$ and a graph $H_i \in
\cP_i$. By 
Lemmas~\ref{ind-her-4} and~\ref{ind-her-5}, the set $\cG_i^{\td}[H_i]
:= \{G \in (\cG_i \cap \bS(\cP_i)) \mid H_i \leq G,\  dec_{\cP_i}(G) = 1 \}$
is also a generating set for $\cP_i$. 

The $*$-join of these $m$ sets is then a generating set for $\cP$,
and we can once again pick out just those graphs that are strict and
have minimum decomposability:

\begin{eqnarray*}
(\cG_1[H_1] * \cdots * \cG_m[H_m])^{\td}  :=  \{G' \in \bS(\cP)\ |\
dec_{\cP}(G') = dec(\cP), \textrm{ and } \forall\, i, \\ 
1\le i\le m,\ \exists\, G_i \in \cG_i^{\td}[H_i],\ G' \in G_1 *
\cdots * G_m\}. 
\end{eqnarray*} 

\stmtnum{Lemma} 
{Let $\cP = \cP_1 \circ \cdots \circ \cP_m$.  Then:
$\cG = (\cG_1[H_1] * \cdots * \cG_m[H_m])^{\td} \subseteq \bS(\cP)$
is a generating set for $\cP$; every
$G \in \cG$ has $dec_{\cP}(G) = dec(\cP)$; and every $G\in\cG$ 
is in the $*$-join
of $m$ $\cP_i$-indecomposable graphs which contain $H_1, \ldots, H_m$
respectively. 
}
\stmtnumend
\newline

We are now ready to prove unique factorisation. As in the hereditary
case, we first show that any two factorisations with exactly
$dec(\cP)$ indecomposable factors must be the same, and then prove
that any factorisation into indecomposable factors must have exactly
$dec(\cP)$ terms. 

\resnum{Theorem\label{ind-her-n}}
{An additive induced-hereditary property $\cP$ can have only one
factorisation with exactly $dec(\cP)$ indecomposable factors. 
}
Let $\cP_1 \circ \cdots \circ \cP_n = \cQ_1 \circ \cdots \circ 
\cQ_n$ be two factorisations of $\cP$ into $n = dec(\cP)$
indecomposable factors. 
%We first order the $\cP_i$'s so that $\cP_i
%\supsetneqq \cP_j \Rightarrow i \geq j$; then for every $i > j$, if
%$\cP_i \not= \cP_j$ there is a graph $X_{i,j} \in \cP_i \setminus
%\cP_j$; if $\cP_i = \cP_j$ define $X_{i,j} := \emptyset$. Let
%$H_{i,0} := \bigcup_{i > j} X_{i,j} \in \cP_i$. 
Label the $\cP_i$'s inductively, beginning with $i=n$, so that, for each
$i$, $\cP_i$ is inclusion-wise maximal among
$\cP_1,\cP_2,\dots,\cP_i$.  For each $i,j$ such that $i>j$, if
$\cP_i\setminus \cP_j\ne\emptyset$, then let $X_{i,j}\in
\cP_i\setminus \cP_j$; if $\cP_i\setminus \cP_j=\emptyset$, then
$\cP_i= \cP_j$ and we set $X_{i,j}$ to be the null graph.  For each
$i$, set $\displaystyle H_{i,0}:=\bigcup_{j<i}X_{i,j}$.  Note
$H_{i,0}\in\cP_i$.  The important point is
 that if $\{L_1,L_2,\dots, L_n\}$ is an
unordered $(\cP_1,\dots,\cP_n)$-partition of some graph 
$G$ such that, for each $i=1,2,\dots,n$, 
$H_{i,0}\leq G[L_i]$, then, by reverse 
induction on $i$ starting at $n$, $G[L_i]\in\cP_i$.

For each $i$, let $\cG_i = \{G_{i,0}, G_{i,1}, G_{i,2}, \ldots \}$ be a
generating set for $\cP_i$. We will construct another generating set
for each $\cP_i$ that will turn out to be contained in some $\cQ_j$;
for graphs $G_{i,s}, H_{i,s}$, we will use the second subscript to
denote which step of our construction we are in. 

For each $s \geq 0$, choose a graph $H'_{s+1} \in (\cG_1[H_{1,s},
G_{1,s}] * \cdots * \cG_n[H_{n,s}, G_{n,s}])^{\td}$, and find an
induced supergraph $H_{s+1}$ whose unique $\cP$-decomposition with
$dec(\cP)$ parts uniformly respects the obvious decomposition of
$H'_{s+1}$. We label as $H_{i,s+1}$ the ind-part of $H_{s+1}$ that
contains the graph from $\cG_i[H_{i,s}, G_{i,s}]$. Then, for each $i$,
$H_{i,0} \leq H_{i,1} \leq H_{i,2} \leq \cdots$ 

For $\cG_i[H_{i,s}, G_{i,s}]$ to be non-empty, we must have $H_{i,s}
\in \cP_i$. We know that the $H_{i,s+1}$'s give an unordered $\{\cP_1,
\ldots, \cP_n\}$-partition of $H_{s+1}$. From the earlier remark, for
$i=1,2,\dots,n$, 
$H_{i,s+1}\in \cP_i$.

The ind-parts of $H_s$ also form its unique $\{\cQ_1, \ldots,
\cQ_n\}$-partition.  Thus, there is some permutation $\varphi_s$ of
$\{1,2,\ldots, n\}$ such that, for each $i$, $H_{i,s} \in
\cQ_{\varphi_s(i)}$. Since there are only finitely many permutations
of $\{1,2,\ldots,n\}$, there must be some permutation $\varphi$ that
appears infinitely often. Now whenever $\varphi_t = \varphi$, we
have $H_{i,1} \leq H_{i,2} \leq \cdots \leq H_{i,t} \in
\cQ_{\varphi(i)}$ so by induced-heredity, for every $s\le t$,
$H_{i,s}$ is in $\cQ_{\varphi(i)}$.  Therefore,  we can take
$\varphi_s = \varphi$, for all $s$. By re-labelling the $\cQ_i$'s, we
can assume $\varphi$ is the identity permutation, so that $H_{i,s} \in
\cQ_i$ for all $i$ and $s$.  

Now for each $i$ and $s$, $G_{i,s-1} \leq H_{i,s}$, so that $\cH_i :=
\{H_{i,1}, H_{i,2}, \ldots \}$ is a generating set for $\cP_i$. But
$\cH_i \subseteq \cQ_i$, so $\cP_i = \la \cH_i\ra \subseteq \cQ_i$. 

By the same reasoning there is a permutation $\tau$ such that $\cQ_i
\subseteq \cP_{\tau(i)}$. We cannot relabel the $\cP_i$'s as well, but
if $\tau^k(i) = i$, then we have $\cP_i \subseteq \cQ_i \subseteq
\cP_{\tau(i)} \subseteq \cQ_{\tau(i)} \subseteq \cP_{\tau^2(i)}
\subseteq \cQ_{\tau^2(i)} \subseteq \cdots \subseteq \cP_{\tau^k(i)} =
\cP_i$, so we must have equality throughout; in particular, $\cP_i =
\cQ_i$ for each $i$. 
\resnumend
\newline

The second piece is analogous to Theorem~\ref{m=n}, but the technical
details are rather different. 

\resnum{Theorem\label{ind-her-m=n}}
{Let 
$\cP_1 \circ \cdots \circ \cP_m$ be a factorisation of the
additive induced-hereditary property $\cP$ into indecomposable additive
induced-hereditary properties. Then $m = dec(\cP)$. 
}
By Lemma~\ref{ind-her-0} any $\cP$-strict graph $G$ has $dec_{\cP}(G) \geq m$, so $dec(\cP) \geq m$. To prove the reverse inequality, we suppose $m < n := dec(\cP)$ and then construct a sequence of graphs until we get a contradiction. When graphs or sets have a double subscript, we will use the second number to denote which step of our construction we are in. For each $i$, we start with some generating set $\cG_i$ consisting only of $\cP_i$-indecomposable $\cP_i$-strict graphs.
\newline

Let $H'_1 \in (\cG_1 * \cdots * \cG_m)^{\td}$, with a corresponding $(\cP_1, \ldots, \cP_m)$-partition  $d_1 = (G'_{1,1}, \ldots, G'_{m,1})$, where each $G'_{i,1}$ is $\cP_i$-strict and $\cP_i$-indecomposable.
By Corollary~\ref{unique-respect} there is a graph $H_1 \geq H'_1$ in
$\bS^{\tD}$ whose ind-parts respect $d_1$ uniformly. That is, denoting
the ind-parts by $H_{1,1}, \ldots, H_{n,1}$, there is a partition
$(J_{1,1},J_{2,1},\dots,J_{m,1})$ of $\{1,2,\ldots,n\}$ such
that $\displaystyle \bigcup_{j \in J_{i,1}} V(H_{j,1})$ induces
$G_{i,1} = t_1 G'_{i,1}$. By additivity of $\cP_i$, $t_1 G'_{i,1}$
is in $\cP_i$, and by Lemma~\ref{ind-her-2} it is $\cP_i$-strict and
$\cP_i$-indecomposable. 

In general suppose we have graphs $H_1, H_2, \ldots, H_{k-1}$ such
that, for each $s=1,2,\dots,k-1$:  
\begin{itemize}
\item[(a)] $H_s$ is $\cP$-strict and uniquely $\cP$-decomposable;
\item[(b)] $dec_{\cP}(H_s) = n$, with ind-parts $H_{1,s}, \ldots, H_{n,s}$;
\item[(c)] $H_1 \leq \cdots \leq H_{k-1}$, with the ind-parts labelled
such that, for $j = 1,\ldots,n$, $H_{j,1} \leq H_{j,2} \leq \cdots
\leq H_{j,k-1}$; 
\item[(d)] there is a partition $(J_{1,s}, J_{2,s},\dots,
J_{m,s})$ of $\{1,2,\ldots,n\}$ such that the union $\displaystyle \bigcup_{j \in
J_{i,s}} V(H_{j,s})$ induces a $\cP_i$-indecomposable graph $G_{i,s}$; and
\item[(e)] for $p < q$, there is at least one $i$ for which
$\displaystyle \bigcup_{j \in 
J_{i,p}} V(H_{j,q})$ does not induce a graph in $\cP_i$.
\end{itemize}

We will find two graphs $H'_k$, $H''_k$ before constructing $H_k$
itself. Because $m < n$, some $G_{i,(k-1)}$ contains more than one
ind-part. Since $G_{i,(k-1)}$ is $\cP_i$-indecomposable, for some
$t$ there is some $H'_k \in t H_{1,(k-1)} * \cdots * t
H_{n,(k-1)}$ for which $\displaystyle \bigcup_{j \in 
J_{i,(k-1)}} t V(H_{j,(k-1)})$ does not induce a graph in
$\cP_i$. Now $H_{k-1} \cup H'_k$ is $\cP$-strict with decomposability $n$
(by Lemma \ref{ind-her-2},
because it contains $H_{k-1}$) and has a $\cP$-decomposition $d'_k$
with $n$ parts, each part being just $(t+1) H_{j,(k-1)}$. By
Corollary \ref{unique-super} we find a
graph $H''_k \geq (H_{k-1} \cup H'_k)$ in $\bS^{\tD}$ whose ind-parts
are just the extension of $d'_k$. 

The graph
$H''_k$ has some $(\cP_1, \ldots, \cP_m)$-partition, and we can extend
this to a partition $d_k$ of $H''_k \cup H'_1$ with $G'_{i,1}$ in the
$\cP_i$-part. We now find a graph $H_k \geq (H''_k \cup H'_1)$ in
$\bS^{\tD}$ whose ind-parts respect this partition. 

Properties (a) and (b) hold for $H_k$ by virtue of being in $\bS^{\tD}$. Since
$H_{k-1} \leq H''_k \leq H_k$, and $H_{k-1}$ is uniquely
$\cP$-decomposable, by Lemma~\ref{ind-her-6} we can label the
ind-parts of $H_k$ to satisfy (c). Condition (e) then follows for any
$p < k-1$, while for $p = k-1$ it holds because of the induced
uniquely $\cP$-decomposable subgraph $H''_k$ (which itself contains
the ``bad'' subgraph $H'_k$). Finally, $d_k$ determines a partition of
$H_k$ with the $i^{th}$ part in $\cP_i$ (because the ind-parts of $H_k$
respect $d_k$) and $\cP_i$-indecomposable (by Lemma~\ref{ind-her-2},
since the $i^{th}$ part contains $G'_{i,1}$). 
\newline

Since there is only a finite number of partitions of
$\{1,2,\ldots,n\}$, at some step $B$ we must end up with a partition
that occurred at some previous step $A < B$. But then (d) contradicts (e). 
\resnumend
%\newline

\stmtnum[\ \cite{uft-2}]{Theorem\label{ind-her-factorn}}
{An additive induced-hereditary property has a factorisation into
$dec(\cP)$ (necessarily indecomposable) additive induced-hereditary
factors. 
} 
\stmtnumend
%\newline

\stmtnum[\ \cite{uft-2}]{Corollary\label{ind-her-irr-ind}}
{An additive induced-hereditary property is irreducible if and only if it is indecomposable.
}
\stmtnumend
%\newline

%+%\stmtnum{Corollary - Unique Factorisation\label{ind-her-uft}}
\stmtnum{Additive Induced-Hereditary Unique Factorisation Theorem\label{ind-her-uft}}
{An additive induced-hereditary property has a unique factorisation
into irreducible additive induced-here\-di\-tary factors, and the number of
factors is exactly $dec(\cP)$. 
}
\stmtnumend
%newline

\section{Related results\label{sec-consequences}}
An important consequence of Theorem~\ref{ind-her-m=n} is that, for irreducible $\cP_i$'s, there are uniquely $(\cP_1, \ldots, \cP_n)$-partitionable graphs, given by Mih\'ok's construction (Corollary~\ref{unique-super}). This was used by Broere and Bucko~\cite{uni-1} to determine when such uniquely partitionable graphs exist if the $\cP_i$'s are allowed to be reducible; and by the first author~\cite{farr-comp} to show that recognising reducible properties is NP-hard, with the exception of the set of bipartite graphs.
\newline

Before proving the uniqueness of the factorisations in~\cite{uft-1} 
and~\cite{uft-2}, we tried without success to prove some related results. 
Their validity for induced-hereditary properties in general
is still open. However, for \emph{additive} induced-hereditary properties 
these results follow quite easily from Unique Factorisation, and we state 
them explicitly below. Note that Proposition~\ref{uft-cor-3} is 
equivalent to Theorem~\ref{ind-her-m=n}.
We also show that unique factorisation for additive hereditary properties
follows both from the result for hereditary compositive properties, and 
from the one for additive induced-hereditary properties.

\stmtnum{Cancellation\label{uft-cor-1}}
{Let $\cA,\cB,\cC$ be additive induced-hereditary properties, $\cA \not=
\cU$. If $\cA \circ \cB = \cA \circ \cC$, then $\cB = \cC$. 
 }
%If $\cA \circ \cB = \cU$ then we must have $\cB = \cU$ (if not, $dec(\cA
%\circ \cB) = dec(\cA) + dec(\cB) < \infty$), and similarly $\cU =
%\cC$. Otherwise the result follows from the unique factorisation of
%$\cA \circ \cB = \cA \circ \cC$. 
\stmtnumend
%\newline

\stmtnum{Corollary\label{uft-cor-2}}
{For additive induced-hereditary properties $\cA', \cA, \cB', \cB$, $\cA \not=
\cU \not= \cB$, if $\cA' \circ \cB' = \cA \circ \cB$, and $\cA'
\subseteq \cA$, $\cB' \subseteq \cB$, then $\cA' = \cA$, $\cB' =
\cB$. 
}
%We have $\cA' \circ \cB' \subseteq \cA' \circ \cB \subseteq \cA \circ
%\cB$, and therefore we must have equality throughout. But then we can
%cancel $\cA'$ in the first equality and $\cB$ in the second (note that
%$\cA' \circ \cB' = \cA \circ \cB \not= \cU$, so $\cB \not= \cU$) to
%get our result. 
\stmtnumend
%\newline

\stmtnum{Proposition\label{uft-cor-3}}
{If $\cQ$ and $\cR$, are additive induced-hereditary properties, then $dec(\cQ \circ\cR) = dec(\cQ) + dec(\cR)$. 
}
%By Theorem~\ref{ind-her-factorn}, $\cQ$ and $\cR$ have factorisations into
%$dec(\cQ)$ and $dec(\cR)$ indecomposable factors respectively; together
%these give a factorisation for $\cQ \circ \cR$, which by
%Theorem~\ref{m=n} must have $dec(\cQ \circ \cR)$ factors.
\stmtnumend
\newline

%+%\stmtnoboxnum{Lemma\label{irr-cop-1}}
%+%{Let $\cQ$ and $\cR$ be two properties, not necessarily induced-hereditary. 
%+%\begin{itemize}
%+%\item[A.] $[\cQ] \circ [\cR] = [\cQ \circ \cR]$. 
%+%\item[B.] $\la \cQ\ra \circ \la \cR\ra = \la \cQ \circ \cR\ra$.$\hfill\square$ 
%+%\end{itemize}
%+%}\stmtnoboxnumend
%+%%\newline

%+%\medskip
A property is {\em strongly irreducible} if it has no factorisation into two non-trivial properties.
We recall that an additive hereditary property is irreducible additive
hereditary 
(respectively, irreducible additive induced-heredi\-tary or irre\-du\-ci\-ble
heredi\-tary 
compositive) if it has no factorisation into two non-trivial additive
hereditary 
(respectively, additive induced-hereditary, or hereditary
compositive) properties.

\risnum{Proposition\label{dc=dec}}
{Let $\cP$ be an additive hereditary property. Then: 
\begin{itemize}
\item[A.] $\cP$ is irreducible additive hereditary 
iff it is strongly irreducible 
%+%iff it is irreducible additive induced-hereditary 
%+%iff it is irreducible hereditary compositive;
\item[B.] $\cP$ has a unique factorisation
into irreducible additive here\-di\-tary factors, and the number of
factors is exactly $dc(\cP) = dec(\cP)$;
\item[C.] if $\cP = \cQ_1 \circ \cdots \circ \cQ_r$, and the 
$\cQ_j$'s are all additive induced-hereditary (or all hereditary compositive),
then they are all additive hereditary.
\end{itemize}
}
%+% A. If $\cP = \cQ \circ \cR$ where $\cQ$ and $\cR$ are additive
%+% hereditary, then these factors are also additive induced-hereditary, 
%+% and hereditary compositive.
%+% Conversely, if $\cP = \cS \circ \cT$, where $\cS$ and $\cT$ are additive
%+% induced-hereditary, then by Lemma~\ref{irr-cop-1}.A, 
%+% $\cP = [\cS] \circ [\cT]$, and $[\cS]$ and $[\cT]$ are additive hereditary.
%+% If $\cS$ and $\cT$ are hereditary compositive, then $dc(\cP) \geq 2$, and
%+% by Theorem 1.1 of~\cite{uft-1} it has a factorisation into $dc(\cP)$
%+% additive hereditary properties.
A. If $\cP = \cS \circ \cT$, where $\cS$ and $\cT$ are any two properties, 
then $\cS+\cT := \{G + H \mid G \in \cS, H \in \cT\}$ is a generating set for $\cP$,
with $dc(\cS+\cT) \geq 2$. By Lemma~\ref{gen-1}, $dc(\cP) \geq 2$, and by Theorem 1.1 of~\cite{uft-1} $\cP$ has a factorisation into $dc(\cP)$ additive hereditary properties.

B. Let $\cP = \cP_1 \circ \cdots \circ \cP_n$, where the $\cP_i$'s are
irreducible additive hereditary properties. Then by A, this must be its
unique factorisation into $dc(\cP)$ irreducible hereditary compositive 
properties, and also its unique factorisation into $dec(\cP)$ irreducible
additive induced-hereditary properties.

C. If we factor each $\cQ_j$ into its irreducible additive induced-hereditary
factors, then by B these irreducible factors are all additive hereditary, so each
$\cQ_j$ is the product of additive hereditary factors.
\risnumend
\newline

An irreducible additive hereditary property is thus 
{\em strongly uniquely factorizable} --- it has exactly one factorisation
even when we allow factors that are not additive or hereditary. 
Szigeti and Tuza~\cite[Problem 4, p.\ 144]{szi-t} 
asked whether this was true for all additive hereditary
properties.  Semani\v{s}in~\cite{sem} gave a class of examples of
additive hereditary properties with non-hereditary factors.
We show in~\cite{uf-her} that the only reducible additive
hereditary property that is strongly uniquely factorisable
is the set of bipartite graphs, which is contained in any reducible
 additive hereditary property. 

In~\cite{uft-1}, however, it is claimed that if the factors of an
additive hereditary property are all
hereditary then they must in fact all be additive hereditary. 
The argument assumes that the factorisation of Theorem~\ref{factorn} is
unique when factoring into any hereditary properties; 
we do not believe that this has been proved --- our proofs of 
uniqueness depend heavily on the additivity of the factors --- so
we leave this as an open question:

\medskip If $\cP=\cQ\circ\cR$, with $\cP$ additive and induced-hereditary,
and $\cQ$ and $\cR$  induced-hereditary, must $\cQ$ and $\cR$ be additive too?
cf.~\cite[Problem 4] {szi-t}
%***\section{Open Problems\label{sec-open}}

%In this section, we present a number of open problems.

%***\begin{enumerate}

%***\item If $\cP=\cQ\circ\cR$, with $\cP$ additive and induced-hereditary,
%***and $\cQ$ and $\cR$  induced-hereditary, must $\cQ$ and $\cR$ be
%%***additive too? 
%***cf.~\cite[Problem 4] {szi-t}

%***\item\label{pnmaximal} 
%***For hereditary $\cP$, when is the join of $n$
%***$\cP$-maximal graphs $\cP^n$-maximal?   In particular, if $G\in
%***{\cM}^*(\cP)$, when is the join of $n$ copies of $G$ in ${\cM}^*(\cP^n)$?

%***\item\label{maximal} Let $\cP_1, \ldots, \cP_n$ be irreducible 
%***additive hereditary properties. When is the join of $n$ indecomposable
%***$\cP_i$-maximal graphs, one from each $P_i$, a 
%***$(P_1 \circ \cdots \circ P_n)$-maximal graph? 
%***\end{enumerate}

%***Questions \ref{pnmaximal} and \ref{maximal} are considered for
%***particular choices of properties in \cite{maxigraph}.  We could
%%***ask, instead, 
%***that the join
%***of the $n$ graphs span \emph{some} $\prod \cP_i$-maximal graph of
%***decomposability $n$. 

\medskip After this paper was first submitted, we discussed this work
with Mih\'ok, who now agrees with our interpretation of the results of
\cite{uft-2} and \cite{uft-1}.  He has also provided a different,
perhaps simpler proof of %+%Theorem \ref{ind-her-m=n}.  
Theorems~\ref{m=n} and \ref{ind-her-m=n}.  
We expect this proof to
appear in some other publication.

\end{document}